# Unconditional existence of conformally hyperbolic Yamabe flows

Mario B. Schulz

*ETH Zürich, Rämistrasse 101, 8092 Zürich, Switzerland*

4 June 2019

We prove global existence of instantaneously complete Yamabe flows on hyperbolic space of arbitrary dimension $m \geq 3$ starting from any smooth, conformally hyperbolic initial metric. We do not require initial completeness or curvature bounds. With the same methods, we show rigidity of hyperbolic space under the Yamabe flow.

Let $(M, g_M)$ be a Riemannian manifold. Let $g_0 = u_0 g_M$ be a conformal metric on $M$ defined by a smooth function $u_0 \colon M \to {]0, \infty[}$. A family $(g(t))_{t \in [0, T[}$ of Riemannian metrics on $M$ is called *Yamabe flow* with initial metric $g_0$, if for all $t \in [0, T[$

$$\begin{cases} \frac{\partial}{\partial t} g(t) = -\mathrm{R}_{g(t)} g(t) \\ \quad g(0) = g_0 \end{cases} \tag{1}$$

where $\mathrm{R}_g$ denotes the scalar curvature of the Riemannian manifold $(M, g)$. Richard Hamilton [12] introduced this flow as alternative approach to the Yamabe problem and showed that solutions to equation (1) exist on any compact manifold $(M, g_0)$ without boundary. Since then, a full theory on compact manifolds was developed with major contributions by Chow [7], Ye [23], Schwetlick and Struwe [18] and Brendle [4, 5].

In dimension $\dim(M) = 2$, where the Yamabe flow coincides with the Ricci flow, Gregor Giesen and Peter Topping [19, 8, 20] obtained existence and uniqueness of instantaneously complete solutions to equation (1) on an arbitrary surface $(M, g_0)$. Instantaneous completeness means that the Riemannian manifold $(M, g(t))$ is geodesically complete for all $t > 0$ even if the initial surface $(M, g_0)$ is incomplete. In [16], the author studied the question whether Giesen and Topping's results generalise to noncompact manifolds of higher dimension and obtained affirmative results for manifolds conformally equivalent to hyperbolic space provided that the conformal factor in the initial metric and the initial scalar curvature are both uniformly bounded from above. In the present paper, we now are able to show existence of instantaneously complete Yamabe flows for *any* conformal initial metric on hyperbolic space $(\mathbb{H}, g_\mathbb{H})$ of dimension $m \geq 3$.

**Theorem 1** (existence). *Let $g_0 = u_0 g_\mathbb{H}$ be any conformal Riemannian metric on hyperbolic space $(\mathbb{H}, g_\mathbb{H})$ of dimension $m \geq 3$. Then, there exists an instantaneously complete Yamabe flow $(g(t))_{t \in [0, \infty[}$ on $\mathbb{H}$ satisfying*

(1) $g(0) = g_0$,

(2) $g(t) \geq m(m-1)t\, g_\mathbb{H}$ *for all $t > 0$.*

*Remark.* Recall that hyperbolic space is a noncompact, simply connected Riemannian manifold of constant sectional curvature $-1$ and scalar curvature $\mathrm{R}_{g_\mathbb{H}} = -m(m-1)$. Theorem 1 improves the existence result obtained in [16], which depended on the uniform upper bounds $g_0 \leq C_0 g_\mathbb{H}$ and $\mathrm{R}_{g_0} \leq K$ on the initial metric and its scalar curvature. Here, we are able to drop these assumptions entirely, so the existence of instantaneously complete Yamabe flows on hyperbolic space of any dimension is true with the same level of generality as in dimension two. It is likely that the proof of Theorem 1 can be generalised to apply on any noncompact manifold $(M, g_M)$ with strictly negative scalar curvature $-\kappa_2 \leq \mathrm{R}_{g_M} \leq -\kappa_1 < 0$ in place of hyperbolic space.

However, there exist (geodesically incomplete) initial manifolds $(M, g_0)$ which do *not* allow any instantaneously complete solution to the Yamabe flow. Conformally flat examples $(M, g_0)$ such as the punctured sphere in dimension $m \geq 3$ are given in [16, Theorem 3]. The incompleteness of Yamabe flows on arbitrary punctured manifolds will be analysed in a forthcoming article.

For the more general class of complete, noncompact background manifolds $(M, g_M)$ with nonpositive, bounded scalar curvature and positive Yamabe invariant, the author proves global existence of complete Yamabe flows in [17], provided the initial metric $g_0 = u_0 g_M$ already is complete with $c_1 \leq u_0 \leq c_2$ for some constants $c_1, c_2 > 0$, such that cases like the punctured sphere are excluded. However, no assumptions on the curvature of $g_0$ are required.

If $g_0$ is any Riemannian metric on some noncompact manifold $M$, then existence of a global Yamabe flow on $M$ with initial metric $g_0$ was shown by Yinglian An and Li Ma [15] provided that $(M, g_0)$ is complete, with Ricci curvature bounded from below and with bounded, nonpositive scalar curvature and also by Li Ma [14] under the assumption that $(M, g_0)$ is complete with nonnegative scalar curvature $\mathrm{R}_{g_0}$ which allows a positive solution $w > 0$ of the equation $-\Delta_{g_0} w = \frac{m-2}{4(m-1)} \mathrm{R}_{g_0}$ in $M$.

Bahuaud and Vertman [2, 1] constructed Yamabe flows starting from spaces with incomplete edge singularities which are preserved along the flow. Very recently, Choi,

Daskalopoulos, and King [6] constructed solutions to the Yamabe flow on $\mathbb{R}^m$ which develop a type II singularity in finite time.

In general, solutions to problem (1) on noncompact manifolds $M$ are not unique. An example on $M = \mathbb{H}$ is the flat metric $g_0 = g_{\mathbb{E}}$ which is conformally equivalent to $g_{\mathbb{H}}$ according to the Poincaré ball model. By Theorem 1, there exists an instantaneously complete Yamabe flow $(g(t))_{t \in [0,\infty[}$ on $\mathbb{H}$ with $g(0) = g_{\mathbb{E}}$. However, the constant flow given by $\bar{g}(t) = g_{\mathbb{E}}$ for all $t$ is also a solution to equation (1) because $\mathrm{R}_{g_{\mathbb{E}}} = 0$. In [16] we conjectured that uniqueness holds in the class of instantaneously complete Yamabe flows and obtained a partial result in the class of rotationally symmetric, instantaneously complete flows. While the conjecture is still open in general, the methods used in the proof of Theorem 1 yield the following result without the assumption of symmetry or completeness.

**Theorem 2** (rigidity)**.** *Let $(\mathbb{H}, g_{\mathbb{H}})$ be hyperbolic space of dimension $m \geq 3$. Let $(g(t))_{t \in [0,T[}$ be a Yamabe flow on $\mathbb{H}$ with $g(0) = g_{\mathbb{H}}$. Then, the flow is uniquely given by $g(t) = (m(m-1)t + 1)g_{\mathbb{H}}$ for all $t \in [0, T[$.*

## Proofs of the main results

Let $u_0 \colon \mathbb{H} \to {]0, \infty[}$ be the conformal factor of the given metric $g_0 = u_0 g_{\mathbb{H}}$ on $\mathbb{H}$. We assume that the restriction of $u_0$ to any smooth, bounded domain $\Omega \subset \mathbb{H}$ is in the Hölder space $\mathrm{C}^{2,\alpha}(\Omega)$ for some $0 < \alpha < 1$. Since the Yamabe flow preserves the conformal class, any Yamabe flow $(g(t))_{t \in [0,T[}$ on $\mathbb{H}$ with $g(0) = g_0$ is of the form $g(t) = u(\cdot, t) g_{\mathbb{H}}$ with a conformal factor $u \colon \mathbb{H} \times [0, T[ \to {]0, \infty[}$ satisfying

$$\begin{cases} \frac{\partial}{\partial t} u = -\mathrm{R}_g u, \\ u(\cdot, 0) = u_0. \end{cases} \tag{2}$$

In dimension $m = \dim(\mathbb{H}) \geq 3$ we may introduce the exponent

$$\eta = \frac{m-2}{4}$$

and the function $U = u^\eta$ in order to express the scalar curvature of the conformal metric $g = u g_{\mathbb{H}}$ by

$$\mathrm{R}_g = U^{-\frac{m+2}{m-2}} \left( \mathrm{R}_{g_{\mathbb{H}}} U - 4 \tfrac{m-1}{m-2} \Delta_{g_{\mathbb{H}}} U \right)$$

and to formulate the following equivalent evolution equations for $U$ respectively $u$.

$$\frac{1}{m-1} \frac{\partial U}{\partial t} = \left( m\eta U + \Delta_{g_{\mathbb{H}}} U \right) U^{-\frac{1}{\eta}}, \tag{3}$$

$$\frac{1}{m-1} \frac{\partial u}{\partial t} = m + \frac{\Delta_{g_{\mathbb{H}}} u}{u} + \frac{(m-6)}{4} \frac{|\nabla u|^2_{g_{\mathbb{H}}}}{u^2}. \tag{4}$$

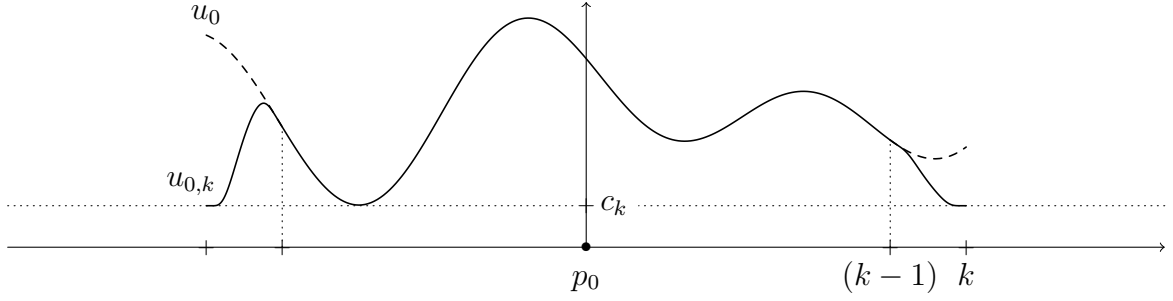

Figure 1: The construction of $u_{0,k}$ given $u_0$.

Let $B_k = B_k(p_0) \subset \mathbb{H}$ be the open metric ball of radius $k$ around some origin $p_0 \in \mathbb{H}$. Let $\chi_k \colon \mathbb{H} \to [0,1]$ be smooth with compact support in $B_k$ satisfying $\chi_k(x) = 1$ for all $x \in B_{k-1}$. For any $k > 2$ we define

$$c_k := \inf_{B_k} u_0 > 0, \tag{5}$$

$$u_{0,k} := (1 - \chi_k) c_k + \chi_k\, u_0, \tag{6}$$

$$\phi_k(t) := c_k + m(m-1)t. \tag{7}$$

Then, $u_{0,k} \in C^{2,\alpha}(B_k)$ coincides with $u_0$ in $B_{k-1}$ and takes the constant value $c_k$ in some neighbourhood of $\partial B_k$ as shown in Figure 1. Moreover, $u_{0,k}$ and $\phi_k$ satisfy the first-order compatibility conditions for equation (4). As shown in [16, Lemma 1.1] there exists some $T_k > 0$ and a solution $u > 0$ of

$$\begin{cases} \dfrac{1}{m-1} \dfrac{\partial u}{\partial t} = m + \dfrac{\Delta_{g_{\mathbb{H}}} u}{u} + \dfrac{(m-6)}{4} \dfrac{|\nabla u|^2_{g_{\mathbb{H}}}}{u^2} & \text{in } B_k \times [0, T_k[, \\ u = \phi_k & \text{on } \partial B_k \times [0, T_k[, \\ u = u_{0,k} & \text{on } B_k \times \{0\}. \end{cases} \tag{8}$$

The following pointwise estimate is analogous to [16, Lemma 1.3].

**Lemma 1.** *Let $u$ be a positive solution to problem (8) with initial and boundary data as given in (6) and (7). Then, for every $0 \leq t < T_k$*

$$\inf_{B_k} u_0 \leq u(\cdot, t) - m(m-1)t \leq \sup_{B_k} u_0.$$

*Proof.* Given any constant $c \in \mathbb{R}$ the function $w(\cdot, t) = u(\cdot, t) - m(m-1)t - c$ satisfies

$$\frac{1}{m-1} \frac{\partial w}{\partial t} - \frac{\Delta_{g_{\mathbb{H}}} w}{u} - \frac{(m-6)\langle \nabla u, \nabla w \rangle_{g_{\mathbb{H}}}}{4u^2} = 0 \quad \text{in } B_k \times [0, T_k[. \tag{9}$$

Since $u > 0$, equation (9) is uniformly parabolic. For $c = \inf_{B_k} u_0$ (respectively $c = \sup_{B_k} u_0$) we have $w \geq 0$ (respectively $w \leq 0$) on $(\partial B_k \times [0, T_k[) \cup (B_k \times \{0\})$ by (7) and the parabolic maximum principle (see [16, Prop. A.2]) implies $w \geq 0$ (respectively $w \leq 0$) in $B_k \times [0, T_k[$. □

**Lemma 2** (global existence on bounded domains). *For every $k > 2$, there exists a unique global solution $0 < u_k \in \mathrm{C}^{2;1}(B_k \times [0, \infty[)$ to problem (8) with $T_k = \infty$, boundary data (7) and initial data (6).*

*Proof.* Since any solution $u \in \mathrm{C}^{2;1}(B_k \times [0, T_k[)$ to problem (8) with $T_k < \infty$ satisfies
$$0 < \inf_{B_k} u_0 \leq u(\cdot, t) \leq \sup_{B_k} u_0 + m(m-1)T_k \quad \text{in } B_k$$
for every $t \in [0, T_k[$ according to Lemma 1, the same approach as in [17, Lemma 1.2] using parabolic DeGiorgi–Nash–Moser theory applies. In fact, a similar argument is used in the proof of Theorem 1. □

The estimate obtained in Lemma 1 is not uniform in $k$ because we do not assume any uniform bounds on $u_0$. To pass to the limit $k \to \infty$ we require local bounds which do not depend on $k$. The nonlinearity of the equation is helpful for upper bounds. Lower bounds however are delicate. We will make use of the following estimate.

**Lemma 3.** *For any real numbers $a, c > 0$ there exists $\lambda > 0$ such that the function $f \colon {]0,1[} \to \mathbb{R}$ given by $f(r) = (1-r^2)^2$ satisfies*
$$f''(r) + \frac{c}{r} f'(r) \geq -\lambda \big(f(r)\big)^{1+a}$$
*for all $r \in {]0,1[}$.*

*Proof.* Since $f'(r) = -4r(1-r^2)$ and $f''(r) = -4 + 12r^2$ we have
$$\begin{aligned}
f''(r) + \frac{c}{r} f'(r) &= (-4 + 12r^2) - 4c(1-r^2) \\
&= 8 - 4(3+c)(1-r^2) \\
&= y\left(\frac{1}{1-r^2}\right)\big(f(r)\big)^{1+a}
\end{aligned}$$
where we introduced the function $y \colon [1, \infty[ \to \mathbb{R}$ given by
$$y(x) = 8x^{2+2a} - 4(3+c)x^{1+2a}.$$
Since the leading term $8x^{2+2p}$ has a positive coefficient, the function $y \colon [1, \infty[ \to \mathbb{R}$ is bounded from below by some constant $-\lambda < 0$ depending only on the parameters $a$ and $c$. □

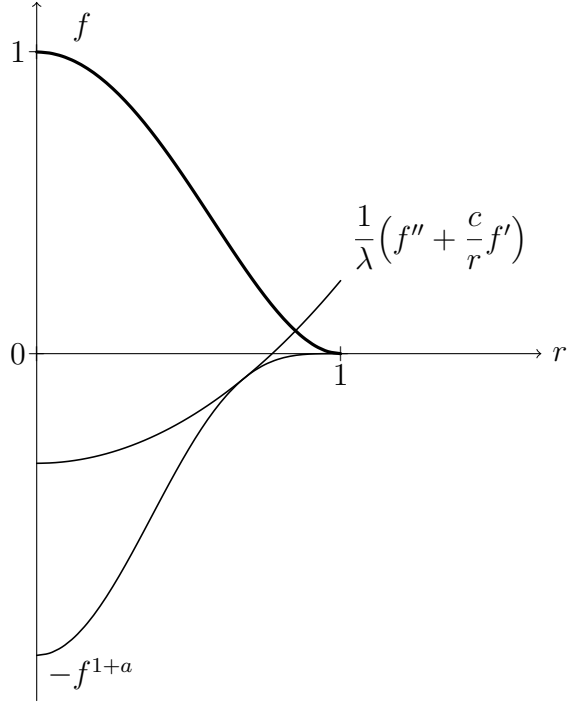

Figure 2: Visualisation of Lemma 3 for $a = 1$, $c = 2$ and $\lambda = 33$.

**Lemma 4** (radial subsolution). *Let $B_1(x_0) \subset \mathbb{H}$ be the open metric unit ball in $(\mathbb{H}, g_\mathbb{H})$ around $x_0 \in \mathbb{H}$ and let $a, b, h_0 > 0$. Then there exists a constant $C > 0$ depending only on $a, b, h_0$ and $m = \dim \mathbb{H}$ such that the map $V \colon B_1 \times [0, t_0[ \to ]0, h_0]$ given by*

$$V(\cdot, t) = \left(h_0^a - Ct\right)^{\frac{1}{a}} \left(1 - r^2\right)^2,$$

*where $t_0 = \frac{1}{C} h_0^a$ and where $r \colon B_1 \to [0, \infty[$ is the Riemannian distance function from $x_0$ in $(\mathbb{H}, g_\mathbb{H})$, satisfies*

$$\tfrac{\partial}{\partial t} V^{1+a} \leq b \Delta_{g_\mathbb{H}} V. \tag{10}$$

*Proof.* A function of the form $V(\cdot, t) = h(t)(f \circ r)$ satisfies (10), if

$$\tfrac{d}{dt} h^{1+a} = -b\lambda h, \qquad -\lambda f^{1+a}(r) \leq \Delta_{g_\mathbb{H}}(f(r))$$

for some $\lambda > 0$. Let $f \colon ]0, 1[ \to \mathbb{R}$ and $\lambda > 0$ be as in Lemma 3 with $c = \frac{m-1}{\tanh(1)}$. Then

$$\Delta_{g_\mathbb{H}}(f(r)) = f''(r) + \frac{m-1}{\tanh(r)} f'(r) \geq f''(r) + \frac{c}{r} f'(r) \geq -\lambda f^{1+a}(r).$$

The equation for $h$ implies $\frac{d}{dt}(h^a) = -\frac{ab\lambda}{a+1}$ which integrates to

$$h^a(t) = h^a(0) - \frac{ab\lambda t}{a+1}.$$

Choosing $h(0) = h_0$ we arrive at

$$V(\cdot, t) = \left(h_0^a - \frac{ab\lambda t}{a+1}\right)^{\frac{1}{a}} \left(1 - r^2\right)^2 \tag{11}$$

which completes the proof with constant $C = \frac{ab\lambda}{a+1}$. $\square$

*Remark.* The equation corresponding to (10) is called *fast diffusion equation* (see [22]) which is well-studied even for domains in Riemannian manifolds of negative curvature: Bonforte, Grillo and Vazquez [3] prove existence of (weak) solutions to the fast diffusion equation in a more general setting and provide more refined estimates of the extinction time. Grillo and Muratori [10] study radial solutions of $\frac{\partial}{\partial t} V^{1+a} = \Delta_{g_\mathbb{H}} V$ on hyperbolic space of dimension $m \geq 3$ in the subcritical range $1 + a < \frac{m+2}{m-2}$ and analyse their fine asymptotics near the extinction time. In the following we will specialise to the critical exponent $1 + a = \frac{m+2}{m-2}$ which corresponds to the Yamabe flow.

It is surprising that the simple profile $f(r) = (1 - r^2)^2$ allows the construction of a compactly supported subsolution to the Yamabe flow on hyperbolic space of any dimension $m \geq 3$. The same approach works on $\mathbb{R}^m$ for $m \geq 3$ if we choose $c = (m-1)$. On manifolds of dimension 2 however, the Yamabe flow behaves differently: According to [9, Theorem A.3] there exist Yamabe flows starting from the flat 2-dimensional unit disc with arbitrarily small extinction time.

**Lemma 5** (local lower bound)**.** *Let $\Omega \subseteq \mathbb{H}$ be any open subset of hyperbolic space of dimension $m \geq 3$ containing the metric ball $B_{r_0} \subset \mathbb{H}$ of radius $r_0 > 1$. Let $(g(t))_{t \in [0,T[}$ be any Yamabe flow on $\Omega$ given by $g(t) = u(\cdot, t) g_\mathbb{H}|_\Omega$. Then, there exists a constant $C_m > 0$ depending only on $m$ and not on $\Omega$ such that for all $t \in [0, T[$*

$$u(\cdot, t) \geq \inf_{B_{r_0}} u(\cdot, 0) - C_m t \quad \text{in } B_{r_0 - 1}.$$

*Proof.* Let $\eta = \frac{m-2}{4}$ as before. According to (3) the function $U = u^\eta$ satisfies

$$\frac{\eta}{(m-1)(\eta+1)} \frac{\partial}{\partial t} U^{1+\frac{1}{\eta}} = m\eta U + \Delta_{g_\mathbb{H}} U \geq \Delta_{g_\mathbb{H}} U$$

in $B_{r_0} \subset \Omega$. Let $x_0 \in B_{r_0-1}$ be arbitrary and let $V \colon B_1(x_0) \times [0, t_0[ \to \mathbb{R}$ be as in Lemma 4 with parameters

$$a = \frac{1}{\eta}, \qquad b = \frac{(m-1)(\eta+1)}{\eta}, \qquad h_0 = \inf_{B_{r_0}} U(\cdot, 0). \tag{12}$$

According to Lemma 2 we may assume $T > t_0$. We consider the difference

$$w := V^{1+\frac{1}{\eta}} - U^{1+\frac{1}{\eta}},$$

define the function $w_+ \colon B_1(x_0) \times [0, t_0[ \to [0, \infty[$ by $w_+(x,t) = \max\{w(x,t), 0\}$ and study the evolution of the quantity

$$J(t) = \int_{B_1(x_0)} w_+(\cdot, t)\, d\mu_{g_\mathbb{H}}.$$

For any $0 < \tau \leq t < t_0$ we have

$$\begin{aligned} J(t) - J(t-\tau) &= \int_{B_1(x_0)} w_+(\cdot,t)\, d\mu_{g_\mathbb{H}} - \int_{B_1(x_0)} w_+(\cdot, t-\tau)\, d\mu_{g_\mathbb{H}} \\ &\leq \int_{\{w(\cdot,t)>0\}} \bigl(w_+(\cdot,t) - w_+(\cdot,t-\tau)\bigr)\, d\mu_{g_\mathbb{H}} \\ &\leq \int_{\{w(\cdot,t)>0\}} \bigl(w(\cdot,t) - w(\cdot,t-\tau)\bigr)\, d\mu_{g_\mathbb{H}}. \end{aligned}$$

Hence,

$$\begin{aligned} \limsup_{\tau \searrow 0} \frac{J(t) - J(t-\tau)}{\tau} &\leq \int_{\{w(\cdot,t)>0\}} \frac{\partial w}{\partial t}(\cdot, t)\, d\mu_{g_\mathbb{H}} \\ &\leq b \int_{\{(V-U)(\cdot,t)>0\}} \Delta_{g_\mathbb{H}}(V-U)(\cdot, t)\, d\mu_{g_\mathbb{H}} \leq 0, \end{aligned} \qquad (13)$$

where we use the following Lemma 6 to obtain the last inequality. We proceed similarly to an argument by Richard Hamilton [11, Lemma 3.1] (see also [16, Lemma A.5]). Let $\varepsilon > 0$ be arbitrary. Estimate (13) implies that there exists $\delta > 0$ such that

$$\forall \tau \in [0, \delta[: \quad J(t) - J(t-\tau) \leq \varepsilon \tau. \qquad (14)$$

We may assume that $\delta \in ]0, t]$ is maximal with this property. By continuity of $t \mapsto J(t)$, estimate (14) extends to

$$J(t) - J(t-\delta) \leq \varepsilon \delta. \qquad (15)$$

If $t - \delta > 0$, we repeat the argument to find $\delta' > 0$ such that

$$\forall \tau \in [0, \delta'[: \quad J(t-\delta) - J(t-\delta-\tau) \leq \varepsilon \tau. \qquad (16)$$

In particular, (15) and (16) can be combined to

$$J(t) - J(t-\delta-\tau) \leq \varepsilon(\delta + \tau)$$

for all $\tau \in [0, \delta'[$ in contradiction to the maximality of $\delta$. Hence, $\delta = t$ and we obtain

$$J(t) - J(0) \leq \varepsilon t.$$

By the choice of $h_0$ we have $J(0) = 0$. Since $\varepsilon > 0$ is arbitrary, $J(t) \leq 0$ follows and implies $U(\cdot,t) \geq V(\cdot,t)$ in $B_1$. In particular, using formula (11) with parameters (12) for $V$, we have

$$U(x_0,t) \geq V(x_0,t) = \left(\inf_{B_{r_0}} u(\cdot,0) - \frac{(m-1)\lambda t}{\eta}\right)^\eta.$$

Since $x_0 \in B_{r_0-1}$ and $t \in {]0,t_0[}$ are arbitrary and $U = u^\eta$, the claim follows with constant $C_m = \frac{m-1}{\eta}\lambda$. □

**Lemma 6.** *Let $\Omega \subset \mathbb{H}$ be a smooth, bounded domain and let $f \in \mathrm{C}^2(\Omega)$ satisfy $f \leq 0$ on $\partial\Omega$. Let $\{f > 0\} := \{x \in \Omega \mid f(x) > 0\}$. Then,*

$$\int_{\{f>0\}} \Delta_{g_\mathbb{H}} f \, d\mu_{g_\mathbb{H}} \leq 0.$$

*Proof.* For any regular value $y \geq 0$ of $f$, the set $\{f > y\} \subset \Omega$ is regular, open and bounded with outer unit normal $\nu$ in the direction of $-\nabla f$. Therefore, we may integrate by parts to obtain

$$\int_{\{f>y\}} \Delta_{g_\mathbb{H}} f \, d\mu_{g_\mathbb{H}} = \int_{\partial\{f>y\}} \langle \nabla f, \nu \rangle_{g_\mathbb{H}} \, d\mu_{g_\mathbb{H}} \leq 0. \tag{17}$$

If $y = 0$ is not a regular value for $f$, we choose a sequence $(y_k)_{k\in\mathbb{N}}$ of regular values for $f$ with $y_k \to 0$ as $k \to \infty$ and pass to the limit in (17). □

The following Lemma about upper bounds is a local version of [16, Proposition 2.1] and complements the local lower bound obtained in Lemma 5. In [16], the estimate is derived from equation (3) for $U = u^\eta$. Here, we give a slightly different proof using equation (4) instead.

**Lemma 7** (local upper bound). *Let $\Omega \subseteq \mathbb{H}$ be any open subset of hyperbolic space of dimension $m \geq 3$ containing the metric ball $B_{r_0} \subset \mathbb{H}$ of radius $r_0 > 1$. Let $(g(t))_{t\in[0,T[}$ be any Yamabe flow on $\Omega$ given by $g(t) = u(\cdot,t)g_\mathbb{H}|_\Omega$. Then, there exists a constant $c_m > 0$ depending only on $m$ and not on $\Omega$ such that for all $t \in [0,T[$*

$$u(\cdot,t) \leq \sup_{B_{r_0}} u(\cdot,0) + (m-1)(m+c_m)t \quad \text{in } B_{r_0-1}.$$

*Proof.* Let $\varphi \colon \Omega \to [0,1]$ be a smooth cutoff function with support in $B_{r_0} \subset \Omega$ such that $\varphi(x) = 1$ for all $x \in B_{r_0-1}$ and such that

$$\left(\frac{(m+2)}{4\varphi}|\nabla\varphi|^2_{g_\mathbb{H}} - \Delta_{g_\mathbb{H}}\varphi\right) \leq c_m \tag{18}$$

in $B_{r_0}$ with some constant $c_m > 0$ depending only on the dimension $m$. Such cutoff functions exist as shown in [16, Lemma A.3–4]. Consider the spatially constant function

$$w(t) = \sup_{B_{r_0}} u(\cdot, 0) + (m-1)(m+c_m)t.$$

Recalling equation (4), but suppressing the index $g_{\mathbb{H}}$ to ease notation of derivatives and inner products, we have

$$\frac{1}{m-1}\frac{\partial}{\partial t}(u\varphi - w)$$
$$= m\varphi - (m+c_m) + \frac{\Delta u}{u}\varphi + \frac{(m-6)}{4}\frac{|\nabla u|^2}{u^2}\varphi$$
$$= m\varphi - (m+c_m) + \frac{\Delta(u\varphi)}{u} - \frac{m+2}{4u}\langle \nabla u, \nabla\varphi\rangle - \Delta\varphi + \frac{m-6}{4u^2}\langle \nabla(\varphi u), \nabla u\rangle$$
$$= m\varphi - (m+c_m) + \frac{\Delta(u\varphi)}{u} - \frac{m+2}{4u\varphi}\langle \nabla(u\varphi), \nabla\varphi\rangle + \frac{m-6}{4u^2}\langle \nabla(\varphi u), \nabla u\rangle$$
$$\quad + \frac{m+2}{4\varphi}|\nabla\varphi|^2 - \Delta\varphi$$
$$\leq \frac{\Delta(u\varphi - w)}{u} - \frac{m+2}{4u\varphi}\langle \nabla(u\varphi - w), \nabla\varphi\rangle + \frac{m-6}{4u^2}\langle \nabla(\varphi u - w), \nabla u\rangle. \quad (19)$$

Since $B_{r_0} \subset \Omega$ is a bounded domain, $u|_{B_{r_0} \times [0,T[}$ is strictly bounded away from zero and from above. Moreover, $u\varphi - w \leq 0$ on $(B_{r_0} \times \{0\}) \cup (\partial B_{r_0} \times [0,T[)$. Hence, the parabolic maximum principle [16, Prop. A.2] applies to inequality (19) and yields $u\varphi - w \leq 0$ in $B_{r_0} \times [0,T[$ which by choice of $\varphi$ implies $u(\cdot, t) \leq w(t)$ in $B_{r_0-1}$ for all $t \in [0,T[$ as claimed. $\square$

*Proof of Theorem 1.* Let $r > 1$ and $T > 1$ be arbitrary but fixed. For every $k \in \mathbb{N}$ let $u_k \colon B_k \times [0, \infty[ \to \,]0, \infty[$ be the solution to problem (8) with boundary data (7) and initial data (6) as given in Lemma 2. Combining the lower bounds from Lemmata 1 and 5, for every $k \geq r+3$ we obtain

$$u_k|_{\overline{B}_{r+2} \times [0,T]} \geq \frac{m(m-1)}{m(m-1) + C_m} \inf_{B_{r+3}} u_0 > 0 \quad (20)$$

where the constant $C_m > 0$ is the same as in Lemma 5. Here we use that we have $\max\{at, b - ct\} \geq \frac{ab}{a+c}$ for any $a, b, t > 0$. In fact, $\max\{at, b - ct\}$ is minimal when $at = b - ct$, that is, when $t = \frac{b}{a+c}$. By Lemma 7, we also have

$$u_k|_{\overline{B}_{r+2} \times [0,T]} \leq \sup_{B_{r+3}} u_0 + (m-1)(m+c_m)T. \quad (21)$$

Recalling $\eta = \frac{m-2}{4}$, we write equation (8) in divergence form

$$\frac{1}{m-1}\frac{\partial u^{\eta+1}}{\partial t} = \frac{m(\eta+1)u^{\eta+1}}{u} + \operatorname{div}_{g_{\mathbb{H}}}\left(\frac{1}{u}\nabla u^{\eta+1}\right) \tag{22}$$

and interpret it as linear parabolic equation for $u_k^{\eta+1}$ in $B_{r+2} \times [0,T]$ with coefficients which are uniformly bounded due to (20) and (21). Since $u_k|_{B_{r+2}\times\{0\}} = u_0|_{B_{r+2}}$ is Hölder continuous by assumption, we may apply parabolic DeGiorgi–Nash–Moser Theory [13, p. 204, Theorem III.10.1] (see also [21, §4]) to equation (22) in order to obtain the interior Hölder bound

$$\|u_k^{\eta+1}\|_{\mathrm{C}^{0,\alpha;0,\frac{\alpha}{2}}(\overline{B}_{r+1}\times[0,T])} \leq C(m,T,u_0|_{B_{r+3}})$$

for some $0 < \alpha < 1$ with a constant $C$ depending only on $m$, $T$, the upper and lower bounds (20) and (21) and the Hölder bound on $u_0|_{B_{r+2}}$, but not on $k$. Together with (20) and (21), it follows that the coefficient $\frac{1}{u_k}$ in equation (22) is Hölder continuous in $\overline{B}_{r+1} \times [0,T]$ satisfying a similar estimate. Since we assume $u_0 \in \mathrm{C}^{2,\alpha}(B_{r+1})$, linear parabolic theory [13, p. 351, Theorem IV.10.1] yields

$$\|u_k^{\eta+1}\|_{\mathrm{C}^{2,\alpha;1,\frac{\alpha}{2}}(\overline{B}_r\times[0,T])} \leq C'(m,T,u_0|_{B_{r+3}}).$$

By compactness of the embedding $\mathrm{C}^{2,\alpha;1,\frac{\alpha}{2}}(\overline{B}_r \times [0,T]) \hookrightarrow \mathrm{C}^{2;1}(\overline{B}_r \times [0,T])$ a subsequence of $\{u_k|_{B_r\times[0,T]}\}_{r+2\leq k\in\mathbb{N}}$ converges to a solution of equation (4) in $\overline{B}_r \times [0,T]$. We repeat this argument to obtain a further subsequence which converges to a solution of equation (4) in $\overline{B}_{2r} \times [0,2T]$.

A diagonal argument allows us to find a subsequence of $\{u_k\}_{k\in\mathbb{N}}$ which converges everywhere to a limit $u \in \mathrm{C}^{2;1}(\mathbb{H} \times [0,\infty[)$ satisfying the Yamabe flow equation (4). Since the uniform lower bound $u(\cdot,t) \geq m(m-1)t$ from Lemma 1 is preserved in the limit, the Yamabe flow given by $g(t) = u(\cdot,t)g_{\mathbb{H}}$ is instantaneously complete. $\square$

*Proof of Theorem 2.* Let $(g(t))_{t\in[0,T[}$ be any Yamabe flow on $\mathbb{H}$ with $g(0) = g_{\mathbb{H}}$. Let $u \colon \mathbb{H} \times [0,T[ \to ]0,\infty[$ be such that $g(t) = u(\cdot,t)g_{\mathbb{H}}$ for all $t \in [0,T[$. In order to obtain a sharp bound on $u$ from below it is convenient to aim for an upper bound on the so-called pressure $v = \frac{1}{u}$ which evolves by the equation (see [16])

$$\tfrac{1}{m-1}\tfrac{\partial}{\partial t}v = -mv^2 + v\Delta_{g_{\mathbb{H}}}v - \tfrac{m+2}{4}|\nabla v|^2_{g_{\mathbb{H}}}. \tag{23}$$

Let $r \colon \mathbb{H} \to [0,\infty[$ denote the Riemannian distance function in $(\mathbb{H}, g_{\mathbb{H}})$ with respect to some origin in $\mathbb{H}$. Given $0 < \varepsilon < \frac{1}{2}$ let $\varphi \colon \mathbb{H} \to ]0,\infty[$ be defined by

$$\varphi = \frac{1}{\cosh(\varepsilon r)}.$$

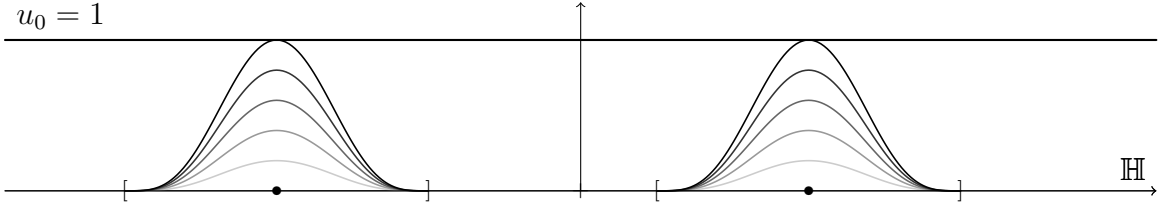

Figure 3: Applying Lemma 5 uniformly.

Then we have $|\nabla\varphi|^2_{g_\mathbb{H}} \leq \varepsilon^2 \varphi^2$ and

$$-\Delta_{g_\mathbb{H}}\varphi = -\frac{\partial^2 \varphi}{\partial r^2} - \frac{(m-1)}{\tanh(r)}\frac{\partial \varphi}{\partial r}$$

$$= \varepsilon^2 \frac{1-\sinh^2(\varepsilon r)}{\cosh^3(\varepsilon r)} + \frac{(m-1)}{\tanh(r)}\frac{\varepsilon \tanh(\varepsilon r)}{\cosh(\varepsilon r)} \leq \varepsilon^2 \varphi + (m-1)\varepsilon\varphi$$

in $\mathbb{H}$. This implies that

$$\begin{aligned}
\tfrac{1}{m-1}\tfrac{\partial}{\partial t}(\varphi v) &= -m\varphi v^2 + v\Delta_{g_\mathbb{H}}(\varphi v) - v^2\Delta_{g_\mathbb{H}}\varphi - 2v\langle\nabla\varphi,\nabla v\rangle_{g_\mathbb{H}} - \tfrac{m+2}{4}\varphi|\nabla v|^2_{g_\mathbb{H}} \\
&\leq -m\varphi v^2 + v\Delta_{g_\mathbb{H}}(\varphi v) + \left(-\Delta_{g_\mathbb{H}}\varphi + \frac{4|\nabla\varphi|^2_{g_\mathbb{H}}}{(m+2)\varphi}\right)v^2 \\
&\leq -m\varphi v^2 + v\Delta_{g_\mathbb{H}}(\varphi v) + m\varepsilon\varphi v^2 \\
&\leq -m(1-\varepsilon)(\varphi v)^2 + v\Delta_{g_\mathbb{H}}(\varphi v)
\end{aligned} \qquad (24)$$

where we used $\varphi \leq 1$ in the last step. Since $u(\cdot, 0) = 1$ we may apply Lemma 5 uniformly in $\mathbb{H}$ and obtain a constant $C_m > 0$ depending only on the dimension $m$ such that $u(\cdot, t) \geq 1 - C_m t$ in $\mathbb{H}$ for all $t \in [0, T[$ as illustrated in Figure 3. This implies that

$$v(\cdot, t) \leq \frac{1}{1 - C_m t}$$

in $\mathbb{H}$ for all $t \in [0, T_0[$, where $T_0 := \min\{T, \frac{1}{C_m}\}$. Hence, the function $(\varphi v)(\cdot, t)$ attains a global maximum in $\mathbb{H}$ and the map $w \colon [0, T_0[ \to\, ]0, \infty[$ given by

$$w(t) = \max_{\mathbb{H}}(\varphi v)(\cdot, t)$$

is well defined. Let $t_0 \in\, ]0, T_0[$ be arbitrary but fixed. Let $q_0 \in \mathbb{H}$ such that

$w(t_0) = (\varphi v)(q_0, t_0)$. By (24), we have

$$\liminf_{\tau \searrow 0} \frac{1}{\tau}\left(\frac{1}{w(t_0)} - \frac{1}{w(t_0 - \tau)}\right) \geq \liminf_{\tau \searrow 0} \frac{1}{\tau}\left(\frac{1}{(\varphi v)(q_0, t_0)} - \frac{1}{(\varphi v)(q_0, t_0 - \tau)}\right)$$

$$= \frac{\partial}{\partial t}\bigg|_{t=t_0} \frac{1}{(\varphi v)(q_0, t)} = \frac{-\frac{\partial}{\partial t}(\varphi v)}{(\varphi v)^2}(q_0, t_0)$$

$$\geq \frac{m-1}{(\varphi v)^2}\Big(m(1-\varepsilon)(\varphi v)^2 - v\Delta_{g_{\mathbb{H}}}(\varphi v)\Big)(q_0, t_0)$$

$$\geq m(m-1)(1-\varepsilon) \tag{25}$$

where we used that $-\Delta_{g_{\mathbb{H}}}(\varphi v)(q_0, t_0) \geq 0$ since $q_0$ is a maximum. As shown in [16, Lemma A.5], estimate (25) implies

$$\frac{1}{w(t)} - \frac{1}{w(0)} \geq m(m-1)(1-\varepsilon)t$$

for every $t \in \,]0, T_0[$ which yields

$$(\varphi v)(\cdot, t) \leq w(t) \leq \frac{1}{m(m-1)(1-\varepsilon)t + 1}$$

since $w(0) = 1$. Letting $\varepsilon \to 0$ and recalling $v = \frac{1}{u}$ we conclude

$$u(\cdot, t) \geq m(m-1)t + 1 \tag{26}$$

for all $t \in [0, T_0[$. By repeating the argument with initial time $T_0$ if necessary, we obtain that estimate (26) holds in fact for all $t \in [0, T[$.

The reverse inequality $u(\cdot, t) \leq m(m-1)t + 1$ is similar to the statement of Lemma 7. In fact, if we choose $\Omega = \mathbb{H}$ and replace the cutoff function $\varphi$ by $\varphi(\varepsilon r)$ in the proof of Lemma 7, then the constant $c_m$ in estimate (18) can be replaced by $\varepsilon c_m$ and we may conclude by letting $\varepsilon \to 0$. □